\documentclass[twoside,11pt,a4paper,leqno]{article}
\usepackage{amsfonts}
\usepackage{amssymb,amsmath,amscd,euscript,verbatim,array}
\usepackage[center,pagestyles]{titlesec}
\usepackage{geometry}
\usepackage{anysize}
\usepackage{fancyhdr}
\usepackage{indentfirst}
\usepackage{graphicx}
\usepackage{color}
\usepackage{ifpdf}
\usepackage{booktabs}
\usepackage{float}
 \marginsize{3.5cm}{3.5cm}{2cm}{2cm}

\titleformat{\section}{\centering\normalsize}{\thesection.}{0.5em}{}
\titleformat{\subsection}{\normalsize\bfseries}{\thesubsection.}{0.5em}{}
\titleformat{\subsubsection}{\normalsize\bfseries}{\thesubsubsection.}{0.5em}{}
\newcommand{\N}{\mathbb{N}}
\newcommand{\Z}{\mathbb{Z}}
\newcommand{\R}{\mathbb{R}}

\newtheorem{Theorem}{Theorem}[section]
\newtheorem{Definition}[Theorem]{Definition}
\newtheorem{Lemma}[Theorem]{Lemma}
\newtheorem{Exercise}[Theorem]{Exercise}
\newtheorem{Proposition}[Theorem]{Proposition}
\newtheorem{Remark}[Theorem]{Remark}

\newcommand{\T}{\mathbb{T}}
\newcommand{\bthm}{\begin{Theorem}}
\newcommand{\ethm}{\end{Theorem}}
\newcommand{\bpr}{\begin{Proposition}}
\newcommand{\epr}{\end{Proposition}}
\newcommand{\blm}{\begin{Lemma}}
\newcommand{\elm}{\end{Lemma}}
\newcommand{\bex}{\begin{Exercise}}
\newcommand{\eex}{\end{Exercise}}
\newcommand{\be}{\begin{equation}}
\newcommand{\ee}{\end{equation}}
\newcommand{\beal}{\begin{aligned}}
\newcommand{\enal}{\end{aligned}}
\newcommand{\brm}{\begin{Remark}}
\newcommand{\erm}{\end{Remark}}
\newcounter{item}[section]

\newcommand{\Proof}{\textbf{Proof}\hspace{0.3cm}}
\newcommand{\End}{\ensuremath{\hfill{\Box}}\\}
\renewcommand{\title}[1]{\begin{center}\textbf{\large #1}\end{center}}
\renewcommand{\author}[1]{\begin{center}\small #1\end{center}}
\renewcommand{\date}[1]{\begin{center}#1\end{center}}

\makeatletter \@addtoreset{equation}{section}
\makeatother
\begin{document}
\vspace{10pt}
\title{DESTRUCTION OF LAGRANGIAN TORUS FOR POSITIVE DEFINITE HAMILTONIAN SYSTEMS}
\vspace{6pt}
\author{\sc Chong-Qing Cheng \& Lin Wang}
\vspace{10pt} \thispagestyle{plain} {\begingroup\makeatletter
\let\@makefnmark\relax
\makeatother\endgroup}
\begin{quote}
\small {\sc Abstract.} For an integrable Hamiltonian
$H_0=\frac{1}{2}\sum_{i=1}^dy_i^2$ $(d\geq 2)$, we show that any
Lagrangian torus with a given unique rotation vector can be
destructed by arbitrarily $C^{2d-\delta}$-small perturbations. In
contrast with it, it has been shown that KAM torus with constant
type frequency persists under $C^{2d+\delta}$-small perturbations
(\cite{P}).
\end{quote}
\begin{quote}
\end{quote}
\begin{quote}
\end{quote} \vspace{25pt}

\section{\sc Introduction and main result}

For exact and area-preserving twist maps on annulus, it was proved
by Herman in \cite{H2} that invariant circles can be destructed by
$C^{3-\delta}$ arbitrarily small perturbations where $\delta$ is a
small positive constant. A variational proof of a similar result was
provided in \cite{W1}. For certain rotation numbers, it was obtained
by Mather (resp. Forni) in \cite{M4} (resp. \cite{Fo}) that the
invariant circles with those rotation numbers can be destroyed by
small perturbations in finer topology respectively. More precisely,
Mather considered Liouville rotation numbers and the topology of the
perturbation induced by $C^\infty$ metric. Forni was concerned about
more special rotation numbers which can be approximated  by rational
ones exponentially and the topology of the perturbation induced by
the supremum norm of $C^\omega$ (real-analytic) function. For
certain positive definite systems with $d$-degrees of freedom,
Herman found that all invariant Lagrangian tori can be destructed by
$C^{d+1-\delta}$ arbitrarily small perturbation in \cite{H3}, where
Lagrangian torus is a natural analogy to invariant circle in
multi-degrees of freedom (see Definition \ref{dd1} below).

In contrast with it, the KAM theory claims the existence of
$d$-dimensional invariant tori in nearly integrable systems with $d$
degrees of freedom. More precisely,  Moser proved that the invariant
circle with constant type frequency of an integrable area-preserving
twist map is persisted under arbitrarily small perturbations in the
$C^{333}$ topology (\cite{Mo1}). Due to the efforts of Moser,
R\"{u}ssman, Herman and P\"{o}schel (\cite{H2,H33,Mo2,Mo3,P,R1,R2}),
it was obtained that certain invariant tori are persisted under
arbitrarily small perturbations in the $C^{2d+\delta}$ topology for
Hamiltonian systems with $d$-degrees of freedom under non-degeneracy
conditions. In particular, Herman proved in \cite{H33} that for
area-preserving twist maps on annulus, certain invariant circles can
be persisted under arbitrarily small perturbations in the $C^3$
topology.

Comparing the results on both directions, it is natural to ask
whether the $C^{2d+\delta}$ condition can be reduced to $C^r$
condition ($r\leq 2d$) to ensure the existence of Lagrangian torus.
In \cite{C}, it was proved that KAM torus with a given rotation
vector does not exist if one carefully construct  perturbations
arbitrarily small in $C^{2d-\delta}$ topology. But this does not
imply non-existence of invariant Lagrangian torus. Indeed, it exists
in the example in \cite{C}.

In this paper we prove the following:
 \begin{Theorem}\label{maint} Given an integrable Hamiltonian $H_0(y)=\frac{1}{2}\sum_{i=1}^dy_i^2$ $(d\geq 2)$,
 a rotation vector $\omega$ and a small positive constant $\delta$,
there exists a sequence of $C^\infty$ Hamiltonians
$\{H_n(x,y)\}_{n\in \N}$ such that $H_n(x,y)\rightarrow H_0(y)$
uniformly in $C^{2d-\delta}$ topology and
 the Hamiltonian flow generated by $H_n(x,y)$ does not admit any
Lagrangian torus with the rotation vector $\omega$, where
$(x,y)\in\T^d\times\R^d$ and $\omega\in\R^d$.
\end{Theorem}

This theorem implies that the rigidity of the Lagrangian torus is as
the same as for the KAM torus. Roughly speaking, if no Lagrangian
torus with the rotation vector $\omega$ survives under an
arbitrarily small perturbation in $C^r$ topology, the maximum of $r$
is closely related to the arithmetic property of the rotation vector
$\omega$. If $\omega$ is a constant type vector, then $r$ is at most
$2d-\delta$. If $\omega$ is a Liouville vector, then $r$ can be
 arbitrarily large. If $\omega$ can be approximated exponentially by
rational vectors, then no Lagrangian torus with the rotation vector
$\omega$ survives under an arbitrarily small perturbation in
$C^\omega$ topology (see \cite{B2}).

In $\text{T}^\ast\T^d$, a submanifold $\mathcal {T}^d$ is called a
Lagrangian torus if it is diffeomorphic to the torus $\T^d$ and the
non-degenerate closed 2-form vanishes on it. For positive definite
Hamiltonian systems, if a Lagrangian torus is invariant under the
Hamiltonian flow, it is then a graph over $\T^d$ (see \cite{BP}). An
example of Lagrangian torus is any KAM torus.
\begin{Definition}\label{dd}$\bar{\mathcal
{T}}^d$ is called a $d$-dimensional KAM torus if
\begin{itemize}
\item $\bar{\mathcal {T}}^d$ is a Lipschitz graph over $\T^d$;
\item $\bar{\mathcal {T}}^d$ is invariant under the Hamiltonian flow
 $\Phi_t^H$ generated by the Hamiltonian function $H$;
\item there exists a diffeomorphism
$\phi:\ \T^d\rightarrow \bar{\mathcal {T}}^d$ such that
$\phi^{-1}\circ\Phi_H^t\circ\phi=R_\omega^t$ for any $t\in \R$,
where $R_\omega^t:\ x\rightarrow x+\omega t$ and $\omega$ is called
the rotation vector of $\bar{\mathcal {T}}^d$.
\end{itemize}
\end{Definition}

Generally, the rotation vector of the Lagrangian torus
 is not well defined if the Lagrangian torus contains several invariant sets with different rotation
vectors. In this paper, we are only concerned with Lagrangian tori
with given unique rotation vectors.
\begin{Definition}\label{dd1}
$\mathcal {T}^d$ is called a  $d$-dimensional Lagrangian torus with
the rotation vector $\omega$ if
\begin{itemize}
\item $\mathcal {T}^d$ is a Lagrangian submanifold;
\item $\mathcal {T}^d$ is invariant for the Hamiltonian
flow $\Phi_H^t$ generated by $H$.
\item each orbit on $\mathcal {T}^d$ has the same rotation vector.
\end{itemize}
\end{Definition}

In order to use a variational method, we recall the definition of
action minimizing orbit for a positive definite Lagrangian system
$L(x,v)$ satisfying superlinear growth with respect to $|v|_x$ (see
\cite{M2}).
\begin{Definition}
An orbit $\gamma\in C^{ac}([t_1,t_2],\T^d)$ (absolutely continuous)
is called an action minimizing curve for $L$ on $[t_1,t_2]$ if
\[\int_{t_1}^{t_2}L(\gamma(t),\dot{\gamma}(t))dt\leq \int_{t_1}^{t_2}L(\bar{\gamma}(t),\dot{\bar{\gamma}}(t))dt,\]
for  every $\bar{\gamma}\in C^{ac}([t_1,t_2],\T^d)$ satisfying
\begin{itemize}
\item $\gamma(t_1)=\bar{\gamma}(t_1)$, $\gamma(t_2)=\bar{\gamma}(t_2)$;
\item $\gamma$ and $\bar{\gamma}$  are in the same homotopy
class, equivalently, their lifts to $\R^d$ connect the same points.
\end{itemize}
\end{Definition}

 In \cite{H4}, it is
proved that each orbit on $\mathcal {T}^d$ is an action minimizing
curve. Actually, in order to prove Theorem \ref{maint}, we only need
the following properties of Lagrangian torus:
\begin{itemize}
\item it is a Lipschitz graph;
\item  every point possesses an action minimizing orbit with rotation vector $\omega$.
\end{itemize}

A rotation vector $\omega\in\R^d$ is called resonant if there exists
$k\in\Z^d$ such that $\langle \omega, k\rangle=0$. Otherwise, it is
non-resonant. The following  approximation of the rotation vector is
found in \cite{C}. For any given vector $\omega\in \R^d$ with $d\geq
2$, there is a sequence of integer vectors $k_n\in\Z^d$ with
$|k_n|\rightarrow\infty$ as $n\rightarrow\infty$ such that
\begin{equation}\label{ap}
\left|\langle\omega, k_n\rangle\right|<\frac{C}{|k_n|^{d-1}},
\end{equation}
where $C$ is a constant independent of $n$ and $|\cdot|$ denotes
Euclidean norm,  i.e. $|k|=(\sum_{j=1}^d k_i^2)^{1/2}$ for $
k=(k_1,k_2,\ldots,k_d)$. (we fix the norm once and for all.) If
$\omega$ is non-resonant, this sequence contains infinitely many
indivisible integer vectors. Here, we call the integer vector
$k\in\Z^d$ indivisible integer vector if $\lambda k\notin\Z^d$
whenever $|\lambda|<1$.

Since a generic hyperbolic perturbation  results in the
transversality between the stable manifold and the unstable
manifold, a Lagrangian torus with resonant rotation vector can be
destructed by analytic perturbation arbitrarily close to zero.
Hence, it is sufficient to consider a Lagrangian torus with
non-resonant rotation vector.

The system we consider consists of one pendulum, a rotator with
$d-1$ degrees of freedom and a perturbation coupling of them. By the
transforation of coordinates (see (\ref{K})), the system with a
given rotation vector $\omega$ can be reduced to the case with the
rotation vector having a small first coordinate and a suitable large
second coordinate. With variational methods, the full action of the
system and the action of the pendulum near the separatrix are
estimated. Moreover, it is obtained that there exists a point such
that no orbit pass through it. Based on the correspondence between
the two cases before and after the transformation (Lemma
\ref{key11}), it is achieved that the original system admits no
Lagrangian torus with the rotation vector $\omega$.

The paper is outlined as follows. In Section 2, we prove Theorem
\ref{maint} for a special case with the rotation vector having a
small first coordinate and a suitable large second coordinate. In
Section 3, the transformation of the coordinates is given and based
on that, we complete the proof of Theorem \ref{maint}.
\section{\sc Destruction of Lagrangian torus with a special rotation
 vector}
After a suitable transformation of coordinates, the problem with a
given rotation vector $\omega$ can be reduced to the case with the
rotation vector having a small first coordinate and a suitable large
second coordinate. In this section, we will prove Theorem
\ref{maint} for that case.

 The Hamiltonian function we consider here is nearly
integrable
\begin{equation}\label{ha}
H_n(q,p)=H_0(p)-P_n(q),\end{equation} where $(q,p)\in
\T^d\times\R^d$. For the sake of simplicity, we
assume\[H_0(p)=\frac{1}{2}|p|^2.\]

Since $H_n$ is strictly convex with respect to $p$, by the Legendre
transformation, the Lagrangian function corresponding to $H_n$ is
\begin{equation}\label{Lo}
L_n(q,\dot{q})=\frac{1}{2}|\dot{q}|^2+P_n(q),
\end{equation}
where $\dot{q}=\frac{\partial H_0}{\partial p}$ and
$q=(q_1,q_2,\ldots,q_d)$.

Let
\[P_n(q)=\frac{1}{n^{a}}(1-\cos q_1)+v_n\left(q_1,q_2\right),\] where
$a$ is a positive constant independent of $n$. For the rotation
vector $\omega=(\omega_1,\ldots,\omega_d)$, $v_n(q_1,q_2)$ is
constructed as follow
 \begin{equation}\label{v}
\begin{cases}
 v_n\ \text{is}\ 2\pi\text{-periodic},\\
\text{supp}\,v_n\cap \{[0,2\pi]\times[-\pi,\pi]\}= B_{R_n}(q^*),\\
\max_{(q_1,q_2)\in [0,2\pi]\times[-\pi,\pi]} v_n=v_n(q^*)={|\omega_1|}^s, \\
{||v_n||}_{C^r}\sim{|\omega_1|}^{s'},
\end{cases}
\end{equation}
where $R_n=\frac{|\omega_1|}{n^2}$, $q^*=(\pi,0)$ and we require
$s'>4$, it can be satisfied for $s$ suitable large. $f\sim g$ means
that $\frac{1}{C}g<f<Cg$ holds for some constant $C>1$,

For (\ref{ha}), we have the following lemma.
\begin{Lemma}\label{key}
For $n$ large enough, the Hamiltonian flow generated by $H_n(q,p)$
does not admit any Lagrangian torus with rotation vector
$\omega=(\omega_1,\ldots,\omega_d)$ satisfying
\[|\omega_1|<n^{-\frac{a}{2}-\epsilon}\quad\text{and}\quad|\omega_2|\sim n,\] where $\epsilon>0$ is independent of $n$.
\end{Lemma}Lemma \ref{key} will be proved with variational method in Subsection 2.3. First
of all, we put it into the Lagrangian formalism. Let
$\sigma_n=n^{-a}$. The Lagrangian function corresponding to
(\ref{Lo}) is
\begin{equation}\label{L1}
\begin{split}
L_n\left(q_1,Q,\dot{q}_1,\dot{Q}\right)=\frac{1}{2}|\dot{Q}|^2+
\frac{1}{2}|\dot{q}_1|^2+\sigma_n(1-\cos(q_1))+v_n(q_1,q_2),
\end{split}\end{equation}
 where $Q=(q_2,\ldots,q_d)$. $L_n(q_1,Q,\dot{q}_1,\dot{Q})$  consists of one pendulum, a rotator with $d-1$ degrees of
freedom and a perturbation coupling of them. The pendulum has the
following Lagrangian
 function
 \begin{equation}\label{simpendul_1}
A_n(q_1,\dot{q}_1)=\frac{1}{2}|\dot{q}_1|^2+\sigma_n(1-\cos(q_1)),
\end{equation}
which corresponds to the Hamiltonian via Legendre transformation
\begin{equation}
h_n(q_1,p_1)=\frac{1}{2}|p_1|^2-\sigma_n(1-\cos q_1).
\end{equation}

\subsection{ The action of the simple pendulum}
Each solution of the Lagrangian equation determined by $A_n$,
denoted by $q_1(t)$, determines an orbit $(q_1(t),p_1(t))$ of the
Hamiltonian flow generated by $h_n$. Each orbit stays in certain
energy level set $(q_1,p_1)\in h_n^{-1}(e)$. Under the boundary
condition that $q_1(t_0)=0,\ q_1(t_1)=\pi$ (or
$q_1(t_1)=\pi,q_1(t_2)=2\pi$), there is a unique correspondence
between $t_1-t_0$ and the energy, denoted by $e(t_1-t_0)$, such that
the determined orbit stays in the energy level set
$h_n^{-1}(e(t_1-t_0))$. More precisely, we have the following lemma.
\begin{Lemma}\label{energ} Let $\bar{q}_1$ be the solution of $A_n$ on
$(t_0,\bar{t}_1)$ satisfying the boundary conditions
\begin{equation*}
\begin{cases}
\bar{q}_1(t_0)=0,\\
\bar{q}_1(\bar{t}_1)=\pi,
\end{cases}
\end{equation*}  $e(\bar{t}_1-t_0)$ be the energy of $\bar{q}_1$, i.e. $(\bar{q}_1,\bar{p}_1)\in h_n^{-1}(e(\bar{t}_1-t_0))$ and
$\omega_1$ be the average speed of $\bar{q}_1$ on $(t_0,\bar{t}_1)$.
If $|\omega_1|<n^{-\frac{a}{2}-\epsilon}$, then
\begin{equation}\label{hh_1}
e(\bar{t}_1-t_0)\sim\sigma_n\exp\left(-\frac{C\sqrt{\sigma_n}}{|\omega_1|}\right),
\end{equation}
where  $\sigma_n=n^{-a}$. \end{Lemma}

 \Proof By the definition,
we have
\[\frac{1}{2}|\dot{\bar{q}}_1|^2-\sigma_n(1-\cos(\bar{q}_1))=e(\bar{t}_1-t_0),\]
hence
\[|\dot{\bar{q}}_1|=\sqrt{2(e(\bar{t}_1-t_0)+\sigma_n(1-\cos(\bar{q}_1)))}.\]
Since $\bar{q}_1(t_0)=0$, then
\[e(\bar{t}_1-t_0)=\frac{1}{2}|\dot{\bar{q}}_1(t_0)|^2\leq \frac{1}{2}|\omega_1|^2\sim n^{-a-2\epsilon},\] which together with $\sigma_n=n^{-a}$ implies that
$\sigma_n/e(\bar{t}_1-t_0)\gg 1$. Since the average speed of
$\bar{q}_1$ is $\omega_1$, by a direct calculation, we have
\[\frac{\pi}{|\omega_1|}=\int_{t_0}^{\bar{t}_1}dt=\int_0^{\pi}\frac{d\bar{q}_1}{\sqrt{2(e(\bar{t}_1-t_0)+\sigma_n(1-\cos(\bar{q}_1)))}}\sim\frac{1}
{\sqrt{\sigma_n}}\ln\left(\frac{\sigma_n}{e(\bar{t}_1-t_0)}\right),\]
moreover,
\[e(\bar{t}_1-t_0)\sim\sigma_n\exp\left(-\frac{C\sqrt{\sigma_n}}{|\omega_1|}\right),\] which completes the proof
of Lemma \ref{energ}.\End

\begin{Remark}
It is easy to see that Lemma \ref{energ} also holds for
\begin{equation*}
\begin{cases}
\bar{q}_1(\bar{t}_1)=\pi,\\
\bar{q}_1(t_2)=2\pi.
\end{cases}
\end{equation*}
Moreover, we denote the difference of time by $\Delta t$, then
 $e(\Delta t)$ is an decreasing function with respect to $\Delta t\geq 0$ for any $\omega_1\in\R$.
\end{Remark}

The following lemma implies that the actions along orbits in the
neighborhood of the separatrix of the  pendulum does not change too
much with respect to a small change in speed (time).

 \begin{Lemma}\label{8_1} Let
$\bar{t}_1, \tilde{t}_1\in [t_0,t_2]$. Let $\bar{q}_1(t)$ be a
solution of $A_n$ on $(t_0,\bar{t}_1)$ and $(\bar{t}_1,t_2)$ with
boundary conditions respectively
 \begin{equation*}\begin{cases}\bar{q}_1(t_0)=0,\\
\bar{q}_1(\bar{t}_1)=\pi,\end{cases}
\begin{cases}\bar{q}_1(\bar{t}_1)=\pi,\\
\bar{q}_1(t_2)=2\pi,\end{cases}\
\end{equation*}
and let $\tilde{q}_1(t)$ be a solution of $A_n$ on
$(t_0,\tilde{t}_1)$ and $(\tilde{t}_1,t_2)$ with boundary conditions
respectively
\begin{equation*}
\begin{cases}\tilde{q}_1(t_0)=0,\\
\tilde{q}_1(\tilde{t}_1)=\pi,\end{cases}
\begin{cases}\tilde{q}_1(\tilde{t}_1)=\pi,\\
\tilde{q}_1(t_2)=2\pi.
\end{cases}
\end{equation*} Let
$\bar{\omega}'_1$ and $\bar{\omega}''_1$ be the average speed of
$\bar{q}_1$ on $(t_0,\bar{t}_1)$ and $(\bar{t}_1,t_2)$ respectively.
Let $\tilde{\omega}'_1$ and $\tilde{\omega}''_1$ be the average
speed of $\tilde{q}_1$ on $(t_0,\tilde{t}_1)$ and
$(\tilde{t}_1,t_2)$ respectively. We set
\[|\omega_1|=\max\left\{|\bar{\omega}'_1|,|\bar{\omega}''_1|,|\tilde{\omega}'_1|,|\tilde{\omega}''_1|\right\}.\] If $|\omega_1|<n^{-\frac{a}{2}-\epsilon}$, then
\begin{equation}\left|\int^{t_2}_{t_0}A_n(\bar{q}_1,\dot{\bar{q}}_1)dt-\int^{t_2}_{t_0}A_n(\tilde{q}_1,\dot{\tilde{q}}_1)dt\right|\leq
C_1|\bar{t}_1-\tilde{t}_1|\sigma_n\exp\left(-\frac{C_2\sqrt{\sigma_n}}{|\omega_1|}\right).\end{equation}
\end{Lemma}

\Proof The proof follows the similar idea of Lemma 4 in \cite{B2}.
Let $q_1(t)$ be a solution of $A_n$ on $(t_0,t_1)$ and $(t_1,t_2)$
with boundary conditions respectively
\begin{equation*}\begin{cases}q_1(t_0)=0,\\
q_1(t_1)=\pi,\end{cases}
\begin{cases}q_1(t_1)=\pi,\\
q_1(t_2)=2\pi.
\end{cases}\end{equation*} We consider the function
\begin{align*}
L(t_1)=&\int_{t_0}^{t_1}A_n(q_1,\dot{q}_1)dt+\int_{t_1}^{t_2}A_n(q_1,\dot{q}_1)dt,\\
=&\int_0^\pi\sqrt{2(e(t_1-t_0)+V(q_1))}dq_1-e(t_1-t_0)(t_1-t_0)\\
&+\int_\pi^{2\pi}\sqrt{2(e(t_2-t_1)+V(q_1))}dq_1-e(t_2-t_1)(t_2-t_1),
\end{align*}
where
\[V(q_1)=\sigma_n(1-\cos(q_1)),\]and $e(\Delta t)$ denotes the energy of the
orbit of the pendulum moving half a turn in time $\Delta t$. The
quantity $e(\Delta t)$ is differentiable with respect to $\Delta t$,
then
\begin{equation}\label{hhh}
\begin{split}
\frac{dL(t_1)}{dt_1}=&\int_0^\pi\frac{\dot{e}(t_1-t_0)}{\sqrt{2(e(t_1-t_0)+V(q_1))}}dq_1-\dot{e}(t_1-t_0)(t_1-t_0)\\
&-e(t_1-t_0)-\int_\pi^{2\pi}\frac{\dot{e}(t_2-t_1)}{\sqrt{2(e(t_2-t_1)+V(q_1))}}dq_1\\
&+\dot{e}(t_2-t_1)(t_2-t_1)+e(t_2-t_1),\\
 =&\int_{t_0}^{t_1}\dot{e}(t_1-t_0)dt-\dot{e}(t_1-t_0)(t_1-t_0)-e(t_1-t_0)\\
 &-\int_{t_1}^{t_2}\dot{e}(t_2-t_1)dt+\dot{e}(t_2-t_1)(t_2-t_1)+e(t_2-t_1),\\
 =&e(t_2-t_1)-e(t_1-t_0).\\
 \end{split}
\end{equation}
Thus, we have \[\left|\frac{dL(t_1)}{dt_1}\right|\leq
|e(t_2-t_1)|+|e(t_1-t_0)|.\]Integrate from $\bar{t}_1$ to
$\tilde{t}_1$ and from (\ref{hh_1}), it follows that
\[\left|\int^{t_2}_{t_0}A_n(\bar{q}_1,\dot{\bar{q}}_1)dt-\int^{t_2}_{t_0}A_n(\tilde{q}_1,\dot{\tilde{q}}_1)dt\right|\leq C_1|\bar{t}_1-\tilde{t}_1|\sigma_n\exp\left(-\frac{C_2\sqrt{\sigma_n}}{|\omega_1|}\right),\]
which completes the proof of Lemma \ref{8_1}. \End
\begin{Remark}\label{aaf}
Since (\ref{hhh}) holds for any $\omega_1\in\R$, combining with the
monotonicity of $e(\Delta t)$, it is easy to draw the following
figure for $L(t)$. $L(t)$ is decreasing for $t\in
(t_0,\frac{t_0+t_2}{2}]$ and increasing for $t\in
[\frac{t_0+t_2}{2},t_2)$. See Fig.1

\input{figure4.TpX}

\end{Remark}

\subsection{ The velocity of the action minimizing orbit} Once the function $q_1(t)$ is fixed, the function $Q(t)$ is the
solution of the Euler-Lagrange equation with the non autonomous
Lagrangian
\begin{equation}\label{LQ}
\hat{L}_n(Q(t),\dot{Q}(t),t)=\frac{1}{2}|\dot{Q}(t)|^2+v_n\left(q_1(t),
q_2(t)\right),
\end{equation}where $Q(t)=(q_2(t),\ldots,q_d(t))$. The first
derivative of $\hat{L}_n(Q(t),\dot{Q}(t),t)$ with respect to $Q$ is
easy to obtained as follows:
\[\frac{\partial \hat{L}_n}{\partial q_2}=\frac{\partial v_n}{\partial q_2}(q_1(t),q_2(t))\quad\text{and}
\quad\frac{\partial \hat{L}_n}{\partial q_i}=0\quad\text{for}\quad
i=3,\ldots,d.\] From the construction of $v_n$ (see (\ref{v})), we
have $||v_n||_{C^r}\sim|\omega_1|^{s'}$. It follows from $s'>4$ that
\[\left|\frac{\partial \hat{L}_n}{\partial q_i}\right|\leq C|\omega_1|^{4}\quad\text{for}\quad
i=2,\ldots,d..\] Based on the periodicity of
$\hat{L}_n(Q(t),\dot{Q}(t),t)$, by Lemma 2 of \cite{BK}, we have the
following estimate.

\begin{Lemma}\label{9_1} Let $(q_1(t),Q(t))$ be the  action minimizing orbit of $L_n$ with
rotation vector $\omega$, then for any $t',t''\in \R$ and $t\in
[t',t'']$ we have
\begin{equation}\label{vv}\left|\dot{Q}(t)-\frac{Q(t'')-Q(t')}{t''-t'}\right|\leq
C|\omega_1|^2.\end{equation}\end{Lemma}

\subsection{ Proof of Lemma \ref{key}}
If the Hamiltonian flow generated by $H_n$ admits a Lagrangian torus
with non-resonant rotation vector $\omega$, then there is a unique
minimal curve $q(t)$ with rotation vector $\omega$ passing through
each $x\in\T^d$ since the Lagrangian torus is a graph. By \cite{H4},
each orbit on the Lagrangian torus is an action minimizing curve.
Hence, it is sufficient to prove the existence of some point in
$\T^d$ where no minimal curve passes through.

Indeed, any minimal curve does not pass through the subspace
$(\pi,0)\times\T^{d-2}$, which implies Lemma \ref{key}. Let
 us assume the contrary, namely, there exists $\bar{t}_1$ such that
\[q_1(\bar{t}_1)=\pi,\quad q_2(\bar{t}_1)=0,\] where
$q(t)=(q_1,q_2,\ldots,q_d)(t)$ is a minimal curve in the universal
covering space $\R^d$. Because of $\omega_1\neq 0$, there exist
$t_0$ and $t_2$ such that
\[q_1(t_0)=0,\quad q_1(t_2)=2\pi.\]
Obviously, $t_0<\bar{t}_1<t_2$.

\noindent\textbf{Claim} $t_2-t_0\sim|\omega_1|^{-1}$.

\noindent\Proof We assume by contradiction that
$t_2-t_0\sim|\bar{\omega}_1|^{-1}$ and without loss of generality,
$\bar{\omega}_1=o(|\omega_1|)$ as $n\rightarrow\infty$. Moreover,
for $[t_0,\bar{t}_1]$ and $[\bar{t}_1,t_2]$, there exists at least
one interval with the length not less than $|\bar{\omega}_1|^{-1}$,
otherwise it is contradicted by $t_2-t_0\sim|\bar{\omega}_1|^{-1}$.
Without loss of generality, one can assume that
$t_2-\bar{t}_1>|\bar{\omega}_1|^{-1}$.

From the definition of rotation vector,
$\lim_{t\rightarrow\infty}(q_1(t)-q_1(\bar{t}_1))/(t-\bar{t}_1)=\omega_1$,
we have that for any $\epsilon>0$, we can find $t_N>t_2$ such that
\begin{equation}
\left|\frac{q_1(t_N)-q_1(\bar{t}_1)}{t_N-\bar{t}_1}-\omega_1\right|\leq\epsilon\quad\text{and}
\quad q_1(t_N)=N\pi,
\end{equation}where $N$ depends on $n$ and $N\gg 2$. Then we have
\[t_N-\bar{t}_1\sim \frac{N\pi}{|\omega_1|}.\]We choose a sequence of times $t_2<\ldots<t_N$ satisfying $q_1(t_i)=i\pi$ for $i\in
\{2,\ldots,N\}$. Moreover, it follows from Pigeon hole principle
that there exists $j\in \{2,\ldots,N-1\}$ such that
\[t_{j+1}-t_j<|\bar{\omega}_1|^{-1},\]where we consider the case
$q_1(t_j)\text{mod}2\pi=0$ and $q_1(t_{j+1})\text{mod}2\pi=\pi$, the
other case is similar. Let $q_1(\bar{t}'_1)=\pi+R_n$ and
$q_1(t'_{j+1})=(j+1)\pi-R_n$, where $R_n$ is the radius of the
support of $v_n$. By the minimality of $q(t)$, there exist positive
constants $C'$ and $C''$ such that
$t_2-\bar{t}'_1>C'|\bar{\omega}_1|^{-1}$ and
$t'_{j+1}-t_j<C''|\bar{\omega}_1|^{-1}$. Hence, we can substitute
$q_1(t)|_{[\bar{t}'_1,t_2]}$ and $q_1(t)|_{[t_j,t'_{j+1}]}$ by the
orbits $\hat{q}_{1'}(t)$, $\hat{q}_{1''}(t)$ of the pendulum $A_n$
(see (\ref{simpendul_1})) with more uniform motion. Correspondingly,
we substitute $q_2(t)|_{[\bar{t}'_1,t_2]}$ and
$q_2(t)|_{[t_j,t'_{j+1}]}$ by the orbits $\hat{q}_{2'}(t)$,
$\hat{q}_{2''}(t)$ of uniform linear motion. From Remark \ref{aaf},
the action of $A_n$ will decrease based on the more uniform motion
of the pendulum. Since the motion of $\hat{q}_{2'}(t)$,
$\hat{q}_{2''}(t)$ are uniform linear , then the actions of
$\frac{1}{2}|\dot{\hat{q}}_{2'}(t)|^2$ and
$\frac{1}{2}|\dot{\hat{q}}_{2''}(t)|^2$ are not greater than the
action of $\frac{1}{2}|\dot{q}_2(t)|^2$. More precisely, similar to
\cite{B2}, we denote
\[T=\frac{1}{2}(t_2-\bar{t}'_1+t'_{j+1}-t_j).\]
Let $s_{j}$ be the closest time to $t'_{j+1}-T$ such that
$q_2(s_j)-q_2(t_j)= l$ where $l\in 2\pi\Z$. Since $|\omega_2|\sim
n$, it follows from Lemma \ref{9_1} that $|t'_{j+1}-T-s_j|\leq C/n$.
Let $s_2=s_j-(t_j-t_2)$. Then $\hat{q}_{1'}(t)$ and
$\hat{q}_{1''}(t)$ are the solutions of $A_n$ on $(\bar{t}'_1,s_2)$
and on $(s_j,t'_{j+1})$ with boundary conditions respectively
\begin{equation*}
\begin{cases}
\hat{q}_{1'}(\bar{t}'_1)=q_1(\bar{t}'_1),\\
\hat{q}_{1'}(s_2)=q_1(t_2)=\pi,
\end{cases}
\quad
\begin{cases}
\hat{q}_{1''}(s_j)=q_1(t_j)=j\pi,\\
\hat{q}_{1''}(t'_{j+1})=q_1(t'_{j+1}).
\end{cases}
\end{equation*}
Correspondingly, we construct $\hat{q}_{2'}(t)$ and
$\hat{q}_{2''}(t)$ respectively
\begin{equation*}
\left\{\begin{array}{ll}\hspace{-0.4em}\hat{q}_{2'}(t)=q_2(\bar{t}'_1)+\frac{q_2(t_2)-q_2(\bar{t}'_1)+l}
{s_2-\bar{t}'_1}(t-\bar{t}'_1),& t\in [\bar{t}'_1,s_2],\\
\hspace{-0.4em}\hat{q}_{2''}(t)=q_2(t_j)+l+\frac{q_2(t'_{j+1})-q_2(t_j)-l}{t'_{j+1}-s_j}(t-s_j),& t\in [s_j,t'_{j+1}].\\
\end{array}\right.
\end{equation*}
 Moreover,
we set
\[\hat{q}_1:=q_1(t)|_{[\bar{t}_1,\bar{t}'_1]}\ast\hat{q}_{1'}(t)|_{[\bar{t}'_1,s_2]}\ast q_1(t-s_2+t_2)|_{[s_2,s_j]}\ast
\hat{q}_{1''}(t)|_{[s_j,t'_{j+1}]}\ast
q_1(t)|_{[t'_{j+1},t_{j+1}]},\]
\[\hat{q}_2:=q_2(t)|_{[\bar{t}_1,\bar{t}'_1]}\ast\hat{q}_{2'}(t)|_{[\bar{t}'_1,s_2]}\ast \left(q_2(t-s_2+t_2)|_{[s_2,s_j]}+l\right)
\ast \hat{q}_{2''}(t)|_{[s_j,t'_{j+1}]}\ast
q_2(t)|_{[t'_{j+1},t_{j+1}]},\] where $\ast$ denotes the
juxtaposition of curves. Let $\hat{Q}=(\hat{q}_2,q_3,\ldots,q_d)$
and $Q=(q_2,q_3,\ldots,q_d)$. Then we have
\[\int_{\bar{t}_1}^{t_N}L_n(\hat{q}_1,\hat{Q},\dot{\hat{q}}_1,\dot{\hat{Q}})dt<\int_{\bar{t}_1}^{t_N}L_n(q_1,Q,\dot{q}_1,\dot{Q})dt,\]
which is contradicted by the minimality of $q(t)$. Hence,
$t_2-t_0\sim|\omega_1|^{-1}$. \End

 Let $\tilde{t}_1$ be the last time before
$\bar{t}_1$ or the first time after $\bar{t}_1$ such that
\[|q_2(\tilde{t}_1)-q_2(\bar{t}_1)|=\pi.\]
Since $|\omega_1|$ is small enough for $n$ large enough and
$|\omega_2|\sim n$, by Lemma \ref{9_1}, we have that for $t\in
[t_0,t_2]$, $|\dot{q}_2(t)|\sim n$,
\[|\tilde{t}_1-\bar{t}_1|\leq \frac{C_0}{n}.\]
 Without loss of generality, one can assume
$\omega_1>0$ and $\omega_2>0$. Consider a solution $\tilde{q}_1$ of
$A_n$ on $(t_0,\tilde{t}_1)$ and on $(\tilde{t}_1,t_2)$ with
boundary conditions respectively
\begin{equation*}
\begin{cases}
\tilde{q}_1(t_0)=q_1(t_0)=0,\\
\tilde{q}_1(\tilde{t}_1)=q_1(\bar{t}_1)=\pi,
\end{cases}
\quad
\begin{cases}
\tilde{q}_1(\tilde{t}_1)=q_1(\bar{t}_1)=\pi,\\
\tilde{q}_1(t_2)=q_1(t_2)=2\pi.
\end{cases}
\end{equation*}
Since $q$ is assumed to be a minimal curve, we have
\begin{equation}\label{star}\int_{t_0}^{t_2}L_n(\tilde{q}_1,Q,\dot{\tilde{q}}_1,\dot{Q})dt-\int_{t_0}^{t_2}L_n(q_1,Q,\dot{q}_1,\dot{Q})dt\geq
0.\end{equation} See Fig.2, where
$x_1=(q_1(\bar{t}_1),q_2(\bar{t}_1))=(\pi,0)$,
$x_0=(q_1(t_0),q_2(t_0))=(0,q_2(t_0))$,
$x_2=(q_1(t_2),q_2(t_2))=(2\pi,q_2(t_2))$, $\tilde{x}'_1=(\pi,-\pi)$
and $\tilde{x}''_1=(\pi,\pi)$.

\input{figure1.TpX}

 $(\tilde{q}_1(t),q_2(t))$ passes through the point $\tilde{x}'_1$ or $\tilde{x}''_1$. Thus, by the construction of $L_n$, we
 obtain from (\ref{star}) that
\begin{equation}\label{AandP_1}
\int_{t_0}^{t_2}A_n(\tilde{q}_1,\dot{\tilde{q}}_1)dt-\int_{t_0}^{t_2}A_n(q_1,\dot{q_1})dt\geq
\int_{t_0}^{t_2}v_n(q_1,q_2)dt-\int_{t_0}^{t_2}v_n(\tilde{q}_1,q_2)dt.
\end{equation}
By the definition of $v_n$ as (\ref{v}), we have the following
claim:

\noindent\textbf{Claim}
$\left(\tilde{q}_1(t),q_2(t)\right)\notin\text{supp}\,v_n$ for $
t\in (t_0,t_2)$.

\noindent\Proof We assume by contradiction that there would exist
$\hat{t}$ such that $(\tilde{q}_1(\hat{t}),q_2(\hat{t}))\in
\text{supp}\,v_n$, without loss of generality, one can assume
$\hat{t}>\tilde{t}_1$. By Lemma \ref{9_1}, for any $t\in
[\tilde{t}_1,\hat{t}]$,
\[\dot{q}_2(t)\leq C_1n,\]hence,
\[\hat{t}-\tilde{t}_1\geq \frac{C_2}{n},\]
where $C_1$, $C_2$ are constants independent of $n$.

Let $\tilde{\omega}'_1$ and $\tilde{\omega}''_1$ be the average
speed of $\tilde{q}_1$ on $(t_0,\tilde{t}_1)$ and
$(\tilde{t}_1,t_2)$ respectively, then
\[\frac{2\pi}{|\omega_1|}=\frac{\pi}{|\tilde{\omega}'_1|}+\frac{\pi}{|\tilde{\omega}''_1|},\]hence,
\[|\tilde{\omega}'_1|\geq \frac{1}{2}|\omega_1|\quad\text{and}\quad|\tilde{\omega}''_1|\geq \frac{1}{2}|\omega_1|,\]
which together with the Euler-Lagrange equation of $\tilde{q}_1(t)$
implies that for any $t\in [\tilde{t}_1,\hat{t}]$,
\[\dot{\tilde{q}}_1(t)\geq C_3|\omega_1|,\]
consequently
\[|\tilde{q}_1(\hat{t})-\tilde{q}_1(\tilde{t}_1)|\geq
C_4\frac{|\omega_1|}{n},\] where $R_n$ is the radius of the support
of $v_n$. Since $R_n=\frac{|\omega_1|}{n^2}$, then we have
\[|\tilde{q}_1(\hat{t})-\tilde{q}_1(\tilde{t}_1)|>R_n,\]
which is contradicted by the assumption
$(\tilde{q}_1(\hat{t}),q_2(\hat{t}))\in \text{supp}\,v_n$.\End

Hence, we have
\[\int_{t_0}^{t_2}v_n(q_1,q_2)dt-\int_{t_0}^{t_2}v_n(\tilde{q}_1,q_2)dt=\int_{t_0}^{t_2}v_n(q_1,q_2)dt.\]
By the construction of $v_n$ and the minimality of $(q_1,Q)$, a
simple calculation shows
\begin{equation}\label{tt}
\int_{t_0}^{t_2}v_n(q_1,q_2)dt\geq {|\omega_1|}^{\lambda},
\end{equation}where $\lambda$ is a positive constant. Consequently, if follows from (\ref{AandP_1}) that
\begin{equation}\label{contr1_1}
\int_{t_0}^{t_2}A_n(\tilde{q}_1,\dot{\tilde{q}}_1)dt-\int_{t_0}^{t_2}A_n(q_1,\dot{q}_1)dt\geq
 {|\omega_1|}^{\lambda}.
\end{equation}

On the other hand, consider a solution $\bar{q}_1$ of $A_n$ on
$(t_0,\bar{t}_1)$ and on $(\bar{t}_1,t_2)$ with boundary conditions
respectively
\begin{equation*}
\begin{cases}
\bar{q}_1(t_0)=q_1(t_0)=0,\\
\bar{q}_1(\bar{t}_1)=q_1(\bar{t}_1)=\pi,
\end{cases}
\quad
\begin{cases}
\bar{q}_1(\bar{t}_1)=q_1(\bar{t}_1)=\pi,\\
\bar{q}_1(t_2)=q_1(t_2)=2\pi.
\end{cases}
\end{equation*}
For $t\in (t_0, \bar{t}_1)$ and $(\bar{t}_1,t_2)$ respectively, the
action of $A_n$ achieves the  minima along $\bar{q}_1(t)$. Thus, we
have
\[\int_{t_0}^{t_2}A_n(q_1,\dot{q}_1)dt\geq\int_{t_0}^{t_2}A_n(\bar{q}_1,\dot{\bar{q}}_1)dt.\]

We compare the action
$\int_{t_0}^{t_2}A_n(\tilde{q}_1,\dot{\tilde{q}}_1)dt$ with the
action $\int_{t_0}^{t_2}A_n(\bar{q}_1,\dot{\bar{q}}_1)dt$  in the
alternative cases, which is  based on the choices of $\tilde{t}_1$.
See Fig.3, where $\bar{t}=\frac{t_0+t_2}{2}$.

\input{figure2.TpX}

\noindent Case 1: $|\bar{t}_1-\bar{t}|\leq \frac{C_0}{n}$.

In this case, the average speed of $\bar{q}_1$ on $(t_0,\bar{t}_1)$
and $(\bar{t}_1,t_2)$ have the same quantity order as $|\omega_1|$.
By $|\tilde{t}_1-\bar{t}_1|\leq \frac{C_0}{n}$, we have
$|\tilde{t}_1-\bar{t}|\leq 2\frac{C_0}{n}$. Hence the average speed
of $\tilde{q}_1$ on $(t_0,\tilde{t}_1)$ and $(\tilde{t}_1,t_2)$ have
also the same quantity order as $|\omega_1|$. Thus, Lemma \ref{8_1}
implies
\[\int_{t_0}^{t_2}A_n(\tilde{q}_1,\dot{\tilde{q}}_1)dt-\int_{t_0}^{t_2}A_n(\bar{q}_1,\dot{\bar{q}}_1)dt\leq
\frac{C_5}{n}\sigma_n\exp\left(-\frac{C_6\sqrt{\sigma_n}}{|\omega_1|}\right).\]

\noindent Case 2: $|\bar{t}_1-\bar{t}|> \frac{C_0}{n}$.

In this case, we take $\tilde{t}_1$ such that
$|\tilde{t}_1-\bar{t}|\leq |\bar{t}_1-\bar{t}|$, which can be
achieved by the suitable choice of the position of
$\tilde{q}_1(\tilde{t}_1)$. More precisely,
\begin{itemize}
  \item  if $\bar{t}_1>\bar{t}+\frac{C_0}{n}$ (Case 2a in Fig.2), we choose $\tilde{t}_1$ as the
last time before $\bar{t}_1$, corresponding to
$(\tilde{q}_1(\tilde{t}_1), q_2(\tilde{t}_1))=\tilde{x}'_1$ in
Fig.2;
  \item  if
$\bar{t}_1<\bar{t}-\frac{C_0}{n}$ (Case 2b in Fig.2), we choose
$\tilde{t}_1$ as the first time after $\bar{t}_1$, corresponding to
$(\tilde{q}_1(\tilde{t}_1), q_2(\tilde{t}_1))=\tilde{x}''_1$ in
Fig.2.
\end{itemize}
 For Case 2a, $\tilde{t}_1\in [\bar{t},\bar{t}_1]$. For Case 2b,  $\tilde{t}_1\in [\bar{t}_1,\bar{t}]$. From Remark \ref{aaf}, it
 follows that $L(\tilde{t}_1)-L(\bar{t})\leq 0$, i.e.
\[\int_{t_0}^{t_2}A_n(\tilde{q}_1,\dot{\tilde{q}}_1)dt-\int_{t_0}^{t_2}A_n(\bar{q}_1,\dot{\bar{q}}_1)dt\leq
0.\] Hence, for any $\bar{t}_1\in (t_0,t_2)$, we can find
$\tilde{t}_1$ such that
\begin{align*}
\int_{t_0}^{t_2}A_n(\tilde{q}_1,\dot{\tilde{q}}_1)dt-\int_{t_0}^{t_2}A_n(q_1,\dot{q}_1)dt&\leq
\int_{t_0}^{t_2}A_n(\tilde{q}_1,\dot{\tilde{q}}_1)dt-\int_{t_0}^{t_2}A_n(\bar{q}_1,\dot{\bar{q}}_1)dt,\\
&\leq
\frac{C_5}{n}\sigma_n\exp\left(-\frac{C_6\sqrt{\sigma_n}}{|\omega_1|}\right).
\end{align*}
Since \[|\omega_1|\leq n^{-\frac{a}{2}-\epsilon}.\] It is easy to
see that for $n$ large enough,
\[\frac{C_5}{n}\sigma_n\exp\left(-\frac{C_6\sqrt{\sigma_n}}{|\omega_1|}\right)\leq
{|\omega_1|}^{\lambda},\]where $\sigma_n=n^{-a}$, which contradicts
to (\ref{contr1_1}) for large $n$. This completes the proof of Lemma
\ref{key}.\End

\section{\sc Destruction of Lagrangian torus with an arbitrary rotation vector}
By (\ref{ap}), for every non resonant rotation vector
$\omega=(\omega_1,\ldots,\omega_d)$ $(d\geq 2)$, there exists a
sequence of integer vector $k_n\in\Z^d$ satisfying $|k_n|\rightarrow
\infty$ as $n\rightarrow \infty$ such that
\[|\langle k_n, \omega\rangle|<\frac{C}{|k_n|^{d-1}}.\]

\begin{Lemma}\label{kkk}
There exists an integer vector $k'_n$ such that
\begin{equation}\label{kn1}
\langle k_n,k'_n\rangle=0,\quad|k'_n|\sim |k_n|\quad\text{and}\quad
|\langle k'_n,\omega\rangle|\sim |k_n|.
\end{equation}
\end{Lemma}
\Proof For $d=2$, Lemma \ref{kkk} holds obviously. In fact, it
suffices to consider $k'_n\in\Z^3$. Let
$k'_n=(k'_{n1},k'_{n2},k'_{n3})$ be the integer vector satisfying
(\ref{kn1}), then for $k'_n\in\Z^d$, one can take
$k'_n=(k'_{n1},k'_{n2},k'_{n3},0,\ldots,0)$ to verify Lemma
\ref{kkk}.

Since $\omega$ is non-resonant, then  $|k_n|\rightarrow\infty$, for
$n\rightarrow \infty$. Let $\Pi$ be the plane orthogonal with
respect to $k_n$. Let $S_{R\alpha}\subset\Pi$ be the sector with
central point $(0,0,0)$, central angle $\alpha$ and radius $R$
satisfying $\alpha\sim 1$. Let the angle between $\omega$ and one of
the radii of $S_{R\alpha}$ be $\beta$, where $\beta\sim 1$ and
$\beta\gg \alpha$.
 Without loss of generality, we
assume $k_{n1}\cdot k_{n2}\cdot k_{n3}\neq 0$. Let
$e_1=(-k_{n2},k_{n1},0)$, $e_2=(0,-k_{n3},k_{n2})$, then $e_1$,
$e_2$ are the generators of the hyperplane $\langle x,k_n\rangle=0$
and satisfy $|e_1|\leq |k_n|$, $|e_2|\leq |k_n|$. Hence, we take
$R\sim |k_n|$, then it is easy to see that the sector $S_{R\alpha}$
contains at least one integer point $m=(m_1,m_2,m_3)$ satisfying
$|m|\sim |k_n|$ for $n$ large enough.

We take $k'_n=(m_1-0,m_2-0,m_3-0)$, then $|k'_n|\sim |k_n|$. Let
$\theta$ be the angle between $k'_n$ and $\omega$, then $\theta\sim
1$. Hence,
\[|\langle
k'_n,\omega\rangle|=|k'_n||\omega||\cos\theta|\sim |k_n|,\]which
together with  $\langle k_n,k'_n\rangle=0$ and $|k'_n|\sim |k_n|$
implies that the integer satisfying (\ref{kn1}) does exist.\End

\subsection{Transformation of coordinates}
By Lemma \ref{kkk}, we choose a sequence of $k_n\in\Z^d$ satisfying
(\ref{ap}) and an integer vector sequence $k'_n$ such that $\langle
k'_n,k_n\rangle=0$, $|k'_n|\sim |k_n|$ and $|\langle
k'_n,\omega\rangle|\sim |k_n|$. In addition, select $d-2$ integer
vectors $l_{n3},\ldots,l_{nd}$ such that
$k_n,k'_n,l_{n3},\ldots,l_{nd}$ are pairwise orthogonal. Let
\begin{equation}\label{K}
K_n=(k_n,k'_n,l_{n3},\ldots,l_{nd})^t. \end{equation} We choose the
transformation of the coordinates
\[q=K_nx.\]
Let $p$ denotes the dual coordinate of $q$ in the sense of  Legendre
transformation, i.e. $p=\frac{\partial L}{\partial\dot{q}}$, it
follows that
\[y=K^t_np,\]where $K^t_n$ denotes the transpose of $K_n$. We
set
\[\Phi_n=\begin{pmatrix}
K_n&\ \\
\ &K^{-t}_n \end{pmatrix},\]then
\[\begin{pmatrix}q\\p\end{pmatrix}=\Phi_n\begin{pmatrix}x\\y\end{pmatrix}.\]
It is easy to verify that \[\Phi_n^tJ_0\Phi_n=J_0,\] where
\[J_0=\begin{pmatrix}
\textbf{0}&\textbf{1}\\
-\textbf{1}&\textbf{0}\end{pmatrix},\] where $\textbf{1}$ denotes a
$d\times d$ unit matrix. Hence, $\Phi_n$ is a symplectic
transformation in the phase space.

\begin{Lemma}\label{key11}
If the Hamiltonian flow generated by $\tilde{H}_n(x,y)$ admits a
Lagrangian torus with rotation vector $\omega$, then the Hamiltonian
flow generated by $H_n(q,p)$ also admits a Lagrangian torus with
rotation vector $K_n\omega$, where $(q,p)^t=\Phi_n(x,y)^t$.
\end{Lemma}

\Proof Let $\tilde{\mathcal {T}}^d$ be the Lagrangian torus admitted
by $\tilde{H}_n(x,y)$, a symplectic form $\Omega$ vanishes on
$\text{T}_x\tilde{\mathcal {T}}^d$ for every $x\in \tilde{\mathcal
{T}}^d$. Since $K_n$ consists of integer vectors, then $\mathcal
{T}^d:=K_n\tilde{\mathcal {T}}^d$ is still a torus. $\Phi_n$ is a
symplectic transformation, hence $\mathcal {T}^d$ is a Lagrangian
torus. From Definition \ref{dd1}, each orbit on $\tilde{\mathcal
{T}}^d$ has the same rotation vector $\omega$. Let
$\tilde{\gamma}(t)$ be a lift of an orbit on $\tilde{\mathcal
{T}}^d$, it follows that
\[\omega=\lim_{t\rightarrow\infty}\frac{\tilde{\gamma}(t)-\tilde{\gamma}(-t)}{2t}.\]
For $\gamma(t)=K_n\tilde{\gamma}(t)$, we have
\begin{align*}
\lim_{t\rightarrow\infty}\frac{\gamma(t)-\gamma(-t)}{2t}&=\lim_{t\rightarrow\infty}\frac{K_n\tilde{\gamma}(t)-K_n\tilde{\gamma}(-t)}{2t};\\
&=K_n\lim_{t\rightarrow\infty}\frac{\tilde{\gamma}(t)-\tilde{\gamma}(-t)}{2t};\\
&=K_n\omega.
\end{align*}
This completes the proof.\End

\subsection{Proof of Theorem \ref{maint}}
We construct $\tilde{H}_n(x,y)$ as follow:
\begin{equation}\label{H}
\tilde{H}_n(x,y)=\frac{1}{2}|y|^2-\tilde{P}_n(x),
\end{equation}
where
\[\tilde{P}_n(x)=\frac{1}{|k_n|^{a+2}}(1-\cos\langle
k_n,x\rangle)+\frac{1}{|k_n|^2}v_n\left(\langle k_n,x\rangle,\langle
k'_n,x\rangle\right),\] where $k'_n$ is the second row of $K_n$ and
$v_n$ is defined by (\ref{v}).  Let $q=K_nx$. In particular, we have
\be\label{cortranfor}
\begin{cases}q_1=\langle k_n,x\rangle,\\
q_2=\langle k'_n,x\rangle.
\end{cases}\ee
By the transformation of coordinates and the Legendre
transformation, the Lagrangian function corresponding to (\ref{H})
is
\begin{equation}\label{ll}
\begin{split}
L_n\left(q_1,Q,\dot{q}_1,\dot{Q}\right)=&\frac{1}{2}\sum^d_{i=3}\frac{|\dot{q}_i|^2}{|l_{ni}|^2}+\frac{|\dot{q}_2|^2}{2|k'_n|^2}\\
&+\frac{1}{|k_n|^2}\left(\frac{1}{2}|\dot{q}_1|^2+\frac{1}{|k_n|^a}(1-\cos(q_1))+v_n(q_1,q_2)\right),
\end{split}
\end{equation} where $Q=(q_2,\ldots,q_d)$.

For the Hamiltonian flow generated by (\ref{H}), by Lemma
\ref{key11}, for  the destruction of the Lagrangian torus
$\tilde{\mathcal {T}}^d$ with rotation vector $\omega$, it is
sufficient to prove that the Euler-Lagrange flow generated by
(\ref{ll}) admits no the Lagrangian torus $\mathcal
{T}^d:=K_n\tilde{\mathcal {T}}^d$ with rotation vector $K_n\omega$.
 Let $K_n\omega=(\omega_1,\ldots,\omega_d)$.

Replacing $n$ by $|k_n|$ in the proof of Lemma \ref{key}, we have
that the Euler-Lagrange flow generated by (\ref{ll}) does not admit
any Lagrangian torus with rotation vector satisfying
\[|\omega_1|<|k_n|^{-\frac{a}{2}-\epsilon}.\]
Here, it should be noted that since $|k'_n|\sim |k_n|$ and the
argument in Section 2 only involves first two components of the
system, then it still holds for (\ref{ll}). From the construction of
$K_n$, $|\omega_1|=|\langle k_n, \omega\rangle|$ which together with
(\ref{ap}) implies
\[|\omega_1|<\frac{C}{|k_n|^{d-1}}.\]
Based on Lemma \ref{key},  it suffices to take
\[\frac{C}{|k_n|^{d-1}}\leq |k_n|^{-\frac{a}{2}-\epsilon},\] which
implies
\begin{equation}\label{a}
a<2d-2-2\epsilon.
\end{equation}

It follows from (\ref{v}) and (\ref{H}) that
\begin{align*}
||\tilde{H}_n&(x,y)-H_0(y)||_{C^r}\\
&=||\tilde{P}_n(x)||_{C^r},\\
&=|k_n|^{-2}\left(|k_n|^{-a}||1-\cos\langle k_n,
x\rangle||_{C^r}+||v_n(\langle k_n,x\rangle, \langle k'_n,
x\rangle)||_{C^r}\right),\\
&\leq |k_n|^{-2}\left(C_1|k_n|^{-a+r}+C_2|k_n|^{-s'(d-1)+r}\right),\\
&\leq C_3\left(|k_n|^{r-a-2}+|k_n|^{r-3(d-1)-2}\right),
\end{align*}
where $\sigma$ is a small positive constant independent of $n$. To
complete the proof of Theorem \ref{maint}, it is enough to make
$r-a-2<0$ and $r-3(d-1)-2<0$, which together with (\ref{a}) implies
\[r<2d-2\epsilon.\]Taking $\delta=3\epsilon$, this completes the proof of
Theorem \ref{maint}.\End

\begin{Remark}
Based on the strategy of the proof, it is easy to see that the
completely integrable part $\frac{1}{2}|y|^2$ of (\ref{H}) can be
generalized to $\frac{1}{2}\langle My,y\rangle$, where
$M=\text{diag}(m_1,\ldots,m_d)$ and $m_i$ $(i=1,\ldots,d)$ are
positive constants independent of $n$. If $M$ is a general positive
definite matrix, some technical difficulties would appear.
\end{Remark}

\vspace{2ex} \noindent\textbf{Acknowledgement} The authors sincerely
thank the referees for their careful reading of the manuscript and
invaluable comments which were very helpful in improving this paper.
The authors also would like to thank Ugo Bessi for illustrations to
his results \cite{B2}. This work is under the support of the NNSF of
China (Grant 11171146), Basic Research Programme of Jiangsu
Province, China (BK2008013), PAPD of Jiangsu Province of China and
Research and Innovation Project for College Graduates of Jiangsu
Province (CXZZ12$\_$0030).

{\sc Department of Mathematics, Nanjing University, Nanjing 210093,
China.}

 {\it E-mail address:} \texttt{chengcq@nju.edu.cn}
\vspace{1em}

{\sc Department of Mathematics, Nanjing University, Nanjing 210093,
China.}

 {\it E-mail address:} \texttt{linwang.math@gmail.com}

\end{document}